\documentclass{llncs}
\usepackage{url}
\usepackage{amssymb}
\usepackage{amsmath}
\usepackage{algorithm}
\usepackage{algorithmic}
\usepackage[left=3cm,top=2.5cm,right=3cm,bottom=2.5cm,nohead]{geometry}
\usepackage{graphicx} 
\usepackage{graphics}
\usepackage{subfigure}
\usepackage{listings}
\DeclareGraphicsExtensions{.pdf}
\usepackage[usenames]{color}

\definecolor{RED}{rgb}{1,0,0}
\definecolor{BLUE}{rgb}{0,0,1}
\begin{document}

\title{Data-Driven Execution of Fast Multipole Methods}
\author{
Hatem Ltaief\inst{1}\and
Rio Yokota\inst{2}
}
\institute{Supercomputing Laboratory\\ \and
Center of Extreme Computing, Division of Mathematical and Computer Sciences and Engineering\\KAUST, Thuwal, KSA\\
\{Hatem.Ltaief,Rio.Yokota\}@kaust.edu.sa
}

\maketitle
\begin{abstract}
Fast multipole methods have $\mathcal{O}(N)$ complexity, are compute bound, 
and require very little synchronization, which makes them a favorable 
algorithm on next-generation supercomputers. Their most common application 
is to accelerate $N$-body problems, but they can also be used to solve 
boundary integral equations. When the particle distribution is irregular 
and the tree structure is adaptive, load-balancing becomes a non-trivial 
question. A common strategy for load-balancing FMMs is to use the work load 
from the previous step as weights to statically repartition the next step. 
The authors discuss in the paper another approach based on data-driven 
execution to efficiently tackle this challenging load-balancing problem.
The core idea consists of breaking the most time-consuming stages
of the FMMs into smaller tasks. The algorithm can then
be represented as a Directed Acyclic Graph (DAG) where
nodes represent tasks, and edges represent dependencies among
them. The execution of the algorithm is performed by asynchronously
scheduling the tasks using the QUARK runtime environment, in a way such that data dependencies
are not violated for numerical correctness purposes. This asynchronous scheduling
results in an out-of-order execution. The performance results of the data-driven
FMM execution outperform the previous strategy and show linear speedup
on a quad-socket quad-core Intel Xeon system.
\end{abstract}
\section{Introduction}
\label{sec:intro}
Bulk synchronous execution/communication models are reaching their limit 
as the amount of concurrency required in high performance 
computing applications increases. This problem is further exacerbated
when the core computational kernels are naturally load imbalanced, which may 
increase the overall idle process time. 

In this paper, we focus on
\emph{fast} $N$-body methods such as the fast multipole method (FMM) on shared-memory
multicore architecture. FMMs have a wide range of applications in astrophysics 
\cite{Dehnen2002}, acoustics \cite{Li2011}, elastodynamics \cite{Chaillat2008}, 
electromagnetics \cite{Vikram2009}, molecular electrostatics \cite{Bajaj2011}, 
and quantum physics \cite{Zhao2007}. Such methods can reduce the complexity of 
the $N$-body problem from $\mathcal{O}(N^2)$ to 
$\mathcal{O}(N)$, while retaining the arithmetic intensity.
On the other hand, they inherently present load balancing issues due to the
irregularity of the data distribution and dynamic nature of the application. FMMs are therefore
a representative class of load imbalanced algorithms for computational science.
One natural solution to this challenging 
problem is the data flow programming model~\cite{dataflow}, which consists of expressing the application
using task-based parallelism. The whole application can then be pictured as 
a directed acyclic graph (DAG), where nodes represent computational tasks
and edges define the data dependencies among them. A dynamic runtime system is then 
employed to efficiently schedule the different tasks over the available
processing units and to ensure the data dependencies are not violated for 
numerical correctness purposes.

Many related works using dynamic scheduling mainly in the dense 
linear algebra area (DLA)~\cite{trefethen,GolubGene}
have been demonstrated in the last few years for solving dense linear systems of
equations~\cite{Agullo2009onesided} as well as eigenvalue and singular value problems
~\cite{Azzam_trd_sc,Azzam_brd_ipdps}. These new \emph{tile} algorithms have been 
integrated into the PLASMA~\cite{PLASMA} high performance DLA library 
targeting x86 multicore platforms, which relies on the dynamic runtime system
QUARK~\cite{quark} to distribute the tasks to worker threads. In the same token, the MAGMA
library~\cite{MAGMA} solves DLA algorithms targeting rather hybrid platforms composed of
x86 multicore enhanced with hardware accelerators (e.g., GPUs) using the StarPU~\cite{Augonnet_2010_ccpe}
framework as the task scheduler. Some efforts
have recently been initiated to extend the PLASMA and MAGMA libraries for distributed 
memory environment, using the dynamic runtime DAGuE~\cite{DAGuE,dplasma}. 
It is noteworthy to mention that there are also other DLA research works based on the same 
fundamental ideas~\cite{CQEQG07,FLAME} i.e., using a data-driven execution framework.
All these numerical frameworks will eventually supersede the state-of-the-art DLA libraries i.e., LAPACK
~\cite{anderson1999} and ScaLAPACK~\cite{ScaLAPACK_1997_guide}, 
for shared and distributed memory systems, respectively.

Similar performance numbers could have been actually 
achieved by using a static scheduler~\cite{Kurzak_2009_ccpe}, 
where each worker thread knows its workload ahead of time. This is mainly because
DLA algorithms are generally well load balanced and compute intensive
enough so that the overhead of the dynamic scheduler is hidden 
by the computational load of the tasks. In fact, the real benefit of using 
the aforementioned dynamic runtime framework for DLA, is to achieve high productivity
in terms of parallel implementation. Therefore, the dynamic feature of the runtime 
is not exploited at all for such well balanced algorithms. It is in this
context that $N$-body methods are an interesting candidate for 
dynamic scheduling of tasks based on data dependencies. $N$-body methods also have computationally 
intensive kernels and would indeed fit well with the data-driven execution model mentioned previously.
The QUARK-enabled FMM implementation presented in this paper achieves a linear
speedup on 16 Intel Xeon cores.

The remainder of this paper is organized as follows: Section
~\ref{sec:fmm} gives a detailed overview of the fast multipole 
method as well as highlighting the load imbalance challenge
and possible ways to fix it.
Section~\ref{sec:runtime} recalls some related works in the dynamic
runtime area, describes the general principles and 
features of the QUARK scheduler and demonstrates how high productivity could be achieved
using this runtime with very small intrusions into the original sequential code. 
Section~\ref{sec:impl} presents the 
implementation details of our FMM code associated with the dynamic scheduler
QUARK. Section~\ref{sec:results} gives the performance results
on a quad-socket quad-core Intel Xeon system (16 cores total). Finally,
Section~\ref{sec:conclusion} summarizes the results of this paper and discusses future work.
\section{Fast Multipole Methods}
\label{sec:fmm}
This Section provides an overview of fast multipole methods, 
describes the inherent load imbalance issue and shows how task scheduling
may efficiently resolve this issue. 

\subsection{Overview}
The fast multipole method (FMM) is a hierarchical $N$-body solver, 
which calculates the interaction of $N$ bodies in $\mathcal{O}(N)$ complexity. 
It has high arithmetic intensity \cite{Yokota2012} and shows good scalability 
on large GPU based systems \cite{Yokota2011}, and also large CPU based systems 
\cite{Rahimian2010}. It possesses a rare combination of linear complexity of 
the algorithm, high arithmetic intensity of the kernels, and locally-dominant 
communication pattern.

\begin{figure}[t]
\centering
\includegraphics[width=0.9\textwidth]{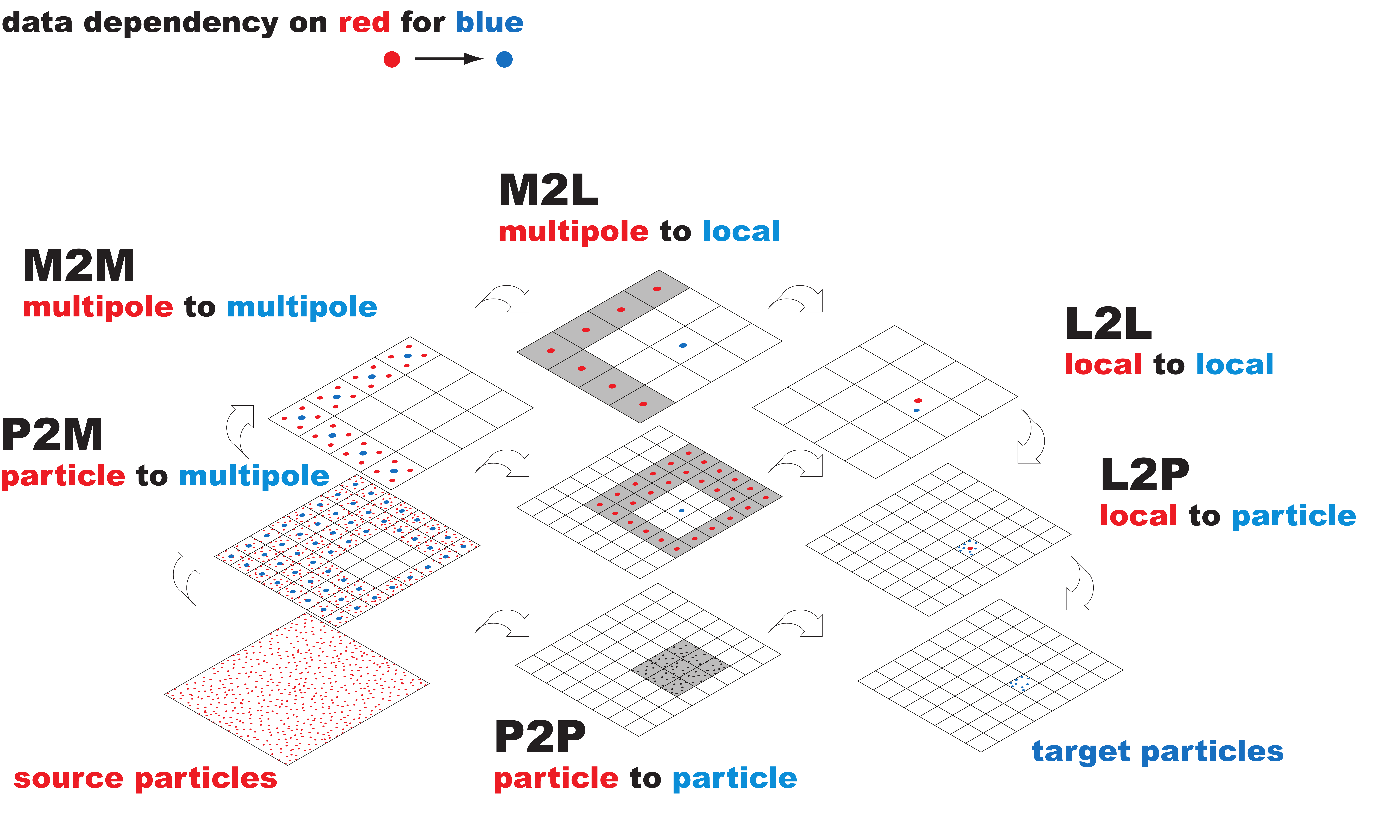}
\caption{Flow of FMM calculation and the individual kernels.}
\label{fig:kernels}
\end{figure}

A schematic of the flow of the FMM algorithm is shown in Figure~\ref{fig:kernels}. 
The FMM consists of six independent stages, some of which have data dependencies on others. 
The domain is partitioned into cells in a hierarchical manner using an octree. 
In this octree structure, we refer to the whole domain as the \textit{root} cell, 
and the smallest cells at the bottom as \textit{leaf} cells. The depth of the tree 
is chosen so that the number of particles per leaf cell remains constant. The particle 
distribution can be irregular in which case the tree structure would become highly adaptive.

The first stage is the particle-to-multipole (P2M) kernel, where the mass/charges 
of the particles are translated into multipole expansions at the center of the leaf cells. 
The next stage is the multipole-to-multipole (M2M) kernel, where the multipole expansion 
at the center of the smaller cells are translated to the center of the larger cells. 
Once the multipole expansions are determined, they can be translated to local expansions 
using the multipole-to-local (M2L) kernel. Note that the M2L kernel can only be performed 
for \textit{well separated} cells as shown in gray in Figure~\ref{fig:kernels}. 
The criteria for choosing well separated cells is based on the ratio between the 
cell size and their distance from each other \cite{Dehnen2002}, or by using the 
parent-neighbor-child relationship in the tree \cite{Cheng1999}. Since the neighboring 
cells at the leaf level will never be handled by the M2L kernel, this part is calculated 
directly by a particle-to-particle (P2P) kernel. Once the local expansions are calculated, 
they are translated to the center of smaller cells by the local-to-local (L2L) kernel. 
Finally, the local expansion at the leaf cell is used to evaluate the solution on the 
particle via the local-to-particle (L2P) kernel. We refer the reader to \cite{Cheng1999} 
for a concise presentation of the mathematical formulae for the individual stages, 
and \cite{Epton1995} for a detailed derivation of these kernels.

Out of the six stages of FMM, the M2L and P2P consume most of the calculation time. 
This is easy to see because P2M and L2P are performed once per leaf cell, and M2M 
and L2L are performed once per cell, while M2L is performed hundreds of times per cell. 
The P2P kernel is also a dominant component of the calculation, since the work load 
between M2L and P2P is always balanced by selecting the number of particles per leaf 
cell to an optimum value. The optimum number of particles per leaf cell depends on the 
efficiency of  the kernel implementation and also the hardware which they run on. 
The P2P kernel has very high arithmetic intensity (Flop/Byte) and is unlikely to 
become bandwidth limited on any architecture, whereas the M2L kernel could become 
bandwidth limited if the order of multipole expansions is small and the architecture 
has low memory bandwidth relative to it's arithmetic capability.

\subsection{Load balancing}
An important difference between the well known domain decomposition problem of meshes, 
and domain decomposition of FMMs, is that each partition in the FMM requires information 
from a global (but hierarchical) halo region, as shown in Figure~\ref{fig:kernels}. 
When the particle distribution is irregular and the tree structure is adaptive, 
load balancing these M2L/P2P kernels becomes a non-trivial task. Different target 
cells will have a different number of source cells to consider, so equally partitioning 
the domain will result in load imbalance. Furthermore, the dynamic nature of the $N$-body 
simulation necessitates frequent repartitioning and load balancing, so the overhead must be small. 
A clever way to solve this dilemma is to record the work load from 
the previous time step and use it to repartition for the present step. 
Such strategies were first used in the early 90's on both sheared memory 
architectures \cite{Singh1993} and distributed memory architectures \cite{Warren1993}. 
The basic idea of using information from the previous step to repartition the present step, 
can be used with any partitioning scheme such as orthogonal recursive 
bisection (ORB) \cite{Warren1992}, partitioning Morton/Hilbert keys \cite{Warren1993}, 
and graph based partitioning \cite{Teng1998}. The work load from the previous step 
can be used as weights in any of these originally unweighted partitioning schemes.

There are many subtleties in the implementation of these partitioning schemes that affect
 their practical usefulness. In order to illustrate these subtleties, we first describe 
the key differences between the two types of tree structures; rectangular binary trees 
resulting from ORB, and cubic octrees. The ORB subdivides the domain into rectangular cuboids, 
while cubic octrees always subdivide into perfect cubes. Since ORB always divides into equal 
number of particles (which could be weighted according to the previous work load) the resulting 
binary tree will be perfectly balanced, and there will be no adaptivity in the tree structure itself. 
When the number of processes is not a power of two, it is a trivial matter to adjust the 
subdivisions so that the tree remains balanced. Therefore, in terms of load balancing, 
this is an ideal tree structure. On the other hand, cubic octrees will result in an adaptive 
tree structure with different depth for different branches if the particle distribution is irregular. 
The advantage of using cubic octrees is the direct correspondence between the proximity of the 
nodes in the tree and the geometrical proximity of the cells which they represent. This means 
Morton/Hilbert keys can be used to determine the geometrical location of the cell and vice versa. 
Therefore, neighbor lists for P2P kernels and well separated lists for M2L kernels can be calculated 
from the Morton/Hilbert keys without traversing the tree and without any explicit information 
about the size and distance of cells. Furthermore, the structured layout of cells permits the 
use of symmetry in the M2L kernels to reduce redundant computation/storage. However, load balancing 
adaptive octrees remains an open area of research \cite{Sundar2008,Cruz2011}. Also, incrementally 
rebalancing the tree is much simpler for the rectangular binary tree since the tree structure 
remains constant and only the cell boundaries are updated, whereas the cubic octree 
requires alteration of the tree structure along with the migration of particles.

In parallel FMMs, the tree structure is used for two separate purposes; a) partitioning and balancing 
the work/communication, and b) traversal that determines the list of cells for M2L/P2P kernels. 
The rectangular binary tree is suitable for the former, while the cubic octree is more appropriate 
for the latter. There are three approaches that can be taken; Case A. 
``Use rectangular binary tree for partitioning and cubic octree for traversal" \cite{Singh1993,Singh1995,Dubinski1996}, 
Case B. ``Use cubic octree for both" \cite{Kurzak2005,Lashuk2009}, and  Case C 
``Use rectangular binary tree for both". There are many techniques 
that can be applied to overcome the weakness of each of these methods, and as far as the 
authors are aware, there is no clear conclusion as to which method is superior.

\subsection{Task scheduling}
Static partitioning of the global tree structure and bulk-synchronous communication of the 
LET are common denominators for the approaches mentioned above. As the amount of concurrency 
required in future architectures continues to increase exponentially, these bulk-synchronous 
execution models will not be able to scale up forever. The convectional technique of using 
work/communication imbalance in the previous step to incrementally rebalance the current step 
can be thought of as a dynamic load balancing scheme in the broad sense. However, the temporal 
granularity of the load balancing is restricted to one time step, and within that time step the 
load balancing is static. One could argue that the particle distribution and tree structure only 
change between time steps, and therefore there is no need to have a finer granularity of load 
balancing if the data structure does not change. This is true if all  threads and all processes 
on all nodes of a large capability machine performed at a constant performance with 0 \% failure rate. 
The fact of the matter is that load balancing on future architectures must also be fault tolerant. 
Furthermore, a single time step can take more than 100 seconds in some 
applications \cite{Yokota2011,Rahimian2010}, and it is definitely worthwhile to load balance at finer
temporal granularity.

Data-driven dynamic task scheduling can solve this issue by stealing work from busy threads, 
while optimizing the data flow at the same time. Since moving data is relatively expensive 
compared to computing data, it is extremely important that the data movement is optimized. In the following section, we will give an overview of QUARK, which is a tool that does precisely
what has been described above.
\section{The QUARK Runtime Environment for Dynamic Task Scheduling}
\label{sec:runtime}
This Section gives an overview of different dynamic runtime schedulers and
highlights then the main features of QUARK, which is the task scheduler
selected to run the FMM code on homogeneous x86 shared-memory multicore architecture.

\subsection{Background}
By now, multicore processors are ubiquitous in both
low-end consumer electronics and high-end servers and supercomputer
large installations. This leads to the emergence of numerous
multithreading frameworks, both open-source and commercial,
embracing the idea of task scheduling: Cilk~\cite{Cilk_manual},
Intel Threading Building Blocks~\cite{intel_tbb,intel_tbb_book}, 
OpenMP (tasking features)~\cite{openmp_30_manual}, just to name a few prominent examples.
From our perspective, one especially important category of
such frameworks are the multithreading systems based on
data flow principles. They represent the computation as a
Direct Acyclic Graph (DAG) and schedule tasks at runtime
through resolution of data hazards: Read after Write (RAW),
Write after Read (WAR) and Write after Write (WAW).
QUeueing And Runtime for Kernels (QUARK) is an example
of such a system. Three others, very similar to some extend, academic
projects are also available: SMPSs~\cite{smpss_manual} from Barcelona Supercomputer
Center, SuperMatrix~\cite{CQEQG07} from the University of Texas at Austin 
and StarPU~\cite{Augonnet_2010_ccpe} from INRIA Bordeaux.

While all four systems have their strength and weaknesses,
QUARK~\cite{quark} has vital extensions for use in a numerical
library as well as at the application level.

\subsection{Principles}

There are many details about the
internals of the scheduler, its dependency analysis, memory
management, and other performance enhancements
that are not covered here. However, information about
an earlier version of this scheduler can be found in~\cite{quark}.

\paragraph{Description of Dependency Types.} 
The authors describe briefly the scheduling techniques of QUARK.
In order for a scheduler to be able to determine dependencies between
the tasks, it needs to know how each task is using
its arguments. Arguments can be VALUE, which are
copied to the task, or they can be INPUT, OUTPUT, or
INOUT, which have the expected meanings. Given the
sequential order that the tasks are added to the scheduler,
and the way that the arguments are used, we can infer
the relationships between the tasks. A task can read a
data item that is written by a previous task (read-after write
RAW dependency); or a task can write a data item
that is written by previous task (write-after-write WAW
dependency); a task can write a data time that is read
by a previous task (write-after-read WAR dependency).
The dependencies between the tasks form an implicit
DAG, however this DAG is never explicitly realized in
the scheduler. The structure is maintained in the way
that tasks are queued on data items, waiting for the
appropriate access to the data.
The tasks are inserted into the scheduler, which stores
them to be executed when all the dependencies are
satisfied. That is, a task is ready to be executed when
all parent tasks have completed. The execution of ready
tasks is handled by worker threads that simply wait
for tasks to become ready and execute them using
a combination of default tasks assignments and work
stealing.

\begin{figure}[bth] 
  \centering
  \includegraphics[trim=0 0 0 0, clip, width=\columnwidth]{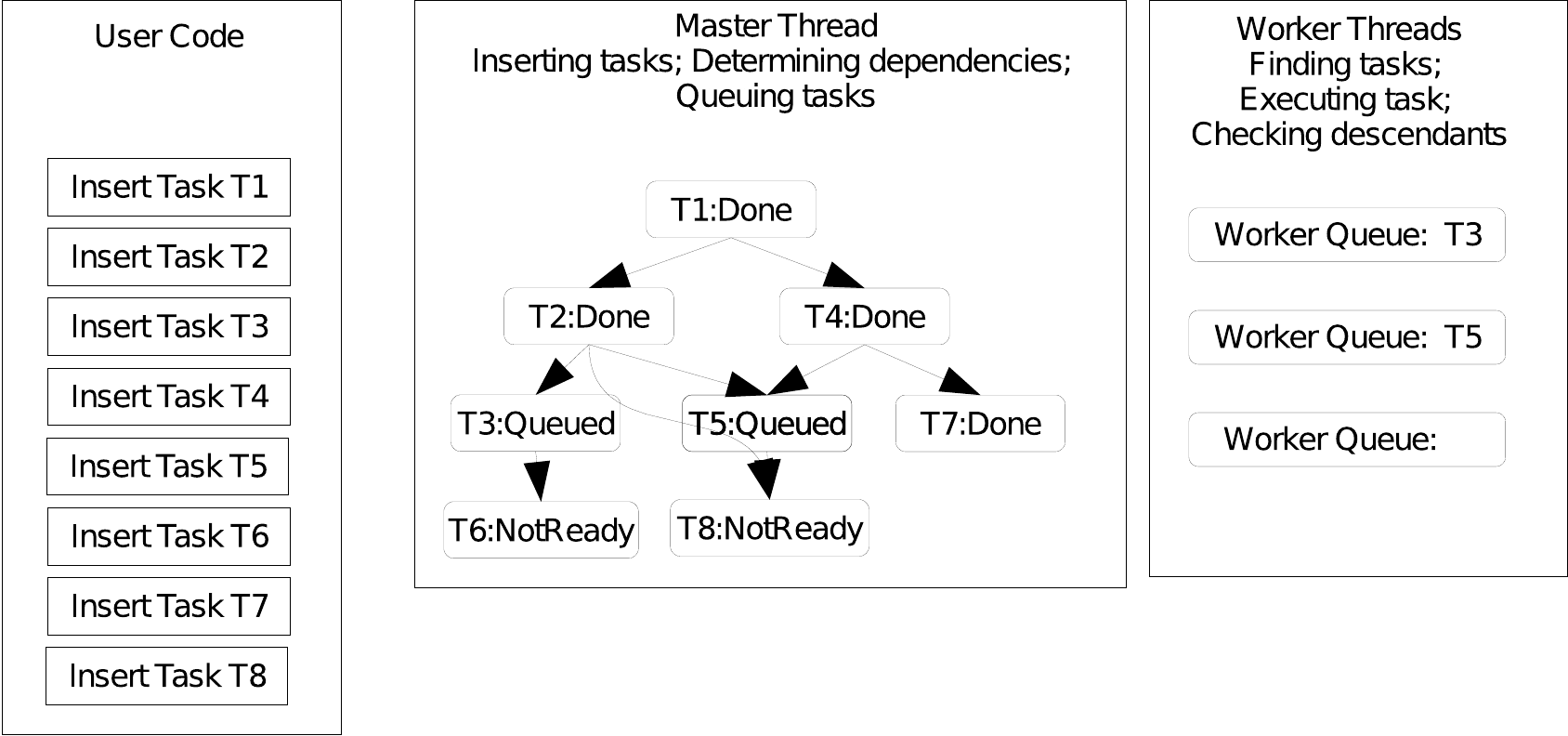}
  \caption{Idealized architecture diagram for the dynamic scheduler.
    Inserted tasks go into a (implicit) DAG based on their
    dependencies.  Tasks can be in NotReady, Queued or Done states.
    Workers execute queued tasks and then determine if any
    descendants have now become ready and can be queued.}
  \label{fig:scheduler-architecture}
\end{figure}

\paragraph{From Sequential Nested-Loop Code to Parallel Execution.}
The scheduler is designed to start from the sequential code (C, C++, Fortran), in which
calls to computational tasks are exposed.
This is intended to make it
easier for algorithm designers to experiment with algorithms and
design new algorithms. 
Each of the calls to the core routines is substituted by a call to
a wrapper that decorates the arguments with their sizes and their
usage (\texttt{INPUT}, \texttt{OUTPUT}, \texttt{INOUT}, \texttt{NODEP},
\texttt{VALUE}). The tasks are inserted into the scheduler, which stores them to be
executed when all the dependencies are satisfied.  That is, a task is
ready to be executed when all parent tasks have completed.  The
execution of ready tasks is handled by worker threads that simply
wait for tasks to become ready and execute them using a combination
of default tasks assignments and work stealing. The thread doing the 
task insertion
is referred to as the master thread.
Under certain circumstances, the master thread will also execute
computational tasks.  Figure~\ref{fig:scheduler-architecture} provides an
idealized overview of the architecture of the dynamic scheduler.

\paragraph{Scheduling a Window of Tasks.}
 
If we were to unfold and retain the entire DAG
of tasks for a large problem, we would be able to perform some
interesting analysis with respect to DAG scheduling and critical
paths.  However, the size of the data structures would quickly grow
overwhelming.  Our solution to this is to maintain a configurable
window of tasks. The implicit DAG is then traversed through
this sliding window, which should be large enough to ensure all 
cores are kept busy. When this window size is reached, the core involved
in inserting tasks does not accept any more tasks until some are
completed. The usage of a window of tasks has implications in how the
loops of an application are unfolded and how much look ahead is
available to the scheduler.

\paragraph{Data Locality and Cache Reuse.}

It has been shown in the past that the reuse of memory caches can lead
to a substantial performance improvement in execution time.  Since we
are working with data structures that should fit in the local caches on
each core, there is a feature which gives the ability to
hint the cache locality behavior.  A parameter in a call 
can be decorated with the LOCALITY flag
in order to tell the scheduler that the data item (parameter) should
be kept in cache if possible.  After a computational core (worker)
executes that task, the scheduler will assign by-default any future
task using that data item to the same core.  Note that the work stealing
can disrupt the by-default assignment of tasks to cores.\\

The next Section discusses further implementation details of the QUARK integration 
into the original FMM source code.
\section{Implementation Details}
\label{sec:impl}
As mentioned in Section \ref{sec:fmm}, the M2L and P2P kernels are the 
dominant part of the FMM calculation. 
Therefore, as a first step, we decided to 
implement the dynamic scheduling by considering as a task each
of these kernels. In our implementation, the M2L and P2P kernels exist inside a 
\textit{dual tree traversal} routine (the \emph{interact} procedure). The 
\textit{dual tree traversal} is a method to find all pairs of well separated 
cells in the octree in $\mathcal{O}(N)$ time \cite{Dehnen2002}. It is more general 
than the commonly used adaptive lists \cite{Cheng1999} because the cells need not 
be perfect cubes. The cells can actually be any shape as long as they are mutually 
exclusive and are grouped hierarchically. Furthermore, the definition of well-separateness 
can be adjusted more smoothly because the definition of neighboring cells is based on 
the distance instead of ``how many cells to skip". This allows the exclusion of the 
corner cells, and yields a list that is closer to a sphere than a cube.

Once a pair of well separated cells are determined, the M2L kernel is called. When the 
dual tree traversal reaches a pair of leaf cells, the P2P kernel is called. In our initial 
implementation, the QUARK interface was placed at the level of these M2L/P2P kernel calls. 
Every time the M2L/P2P kernel was called, this task would be queued by itself. Then, QUARK 
would resolve the data dependency and dynamically schedule the task onto an available 
worker thread. However, it turned out that the overhead of the task scheduler was too large to 
be implemented at this granularity. In other words, there were too many tiny 
tasks to schedule. Note that this is not a matter of arithmetic intensity, but 
the absolute size of both the operations and data being too small.

In order to address the issue above, we implemented a mechanism, which could schedule 
tasks of arbitrary (controllable) granularity. This was done by interfacing QUARK 
at a higher level in the tree traversal. 
Grouping fine-grained tasks naturally follow good data
locality guidelines by retaining data in caches for as long
as the computation continues. QUARK does not migrate
threads that are in the state of execution. And the look-up
and update of the data dependence information takes place
only once per each merged \emph{super} task. The savings in overhead
will come from both lack of cache memory pollution from
QUARK's internal data structures and elimination of tens
if not hundreds of instructions that help QUARK make
scheduling decisions and keep its data in a consistent state.
For example, in the first step of a dual 
tree traversal, the target cell is split into eight child cells and eight new pairs are 
formed between the new target cells and old source cell. If each of these are 
scheduled as tasks on QUARK, the overhead of scheduling will be much smaller than 
the actual task of traversing each of these large subtrees. Thus, there is a tradeoff 
between the small overhead of coarse grain scheduling, and the load balancing that 
fine grain scheduling can offer.

\begin{figure}[tb]
\subfigure[Mutual interaction.]{
\includegraphics[width=0.4\textwidth]{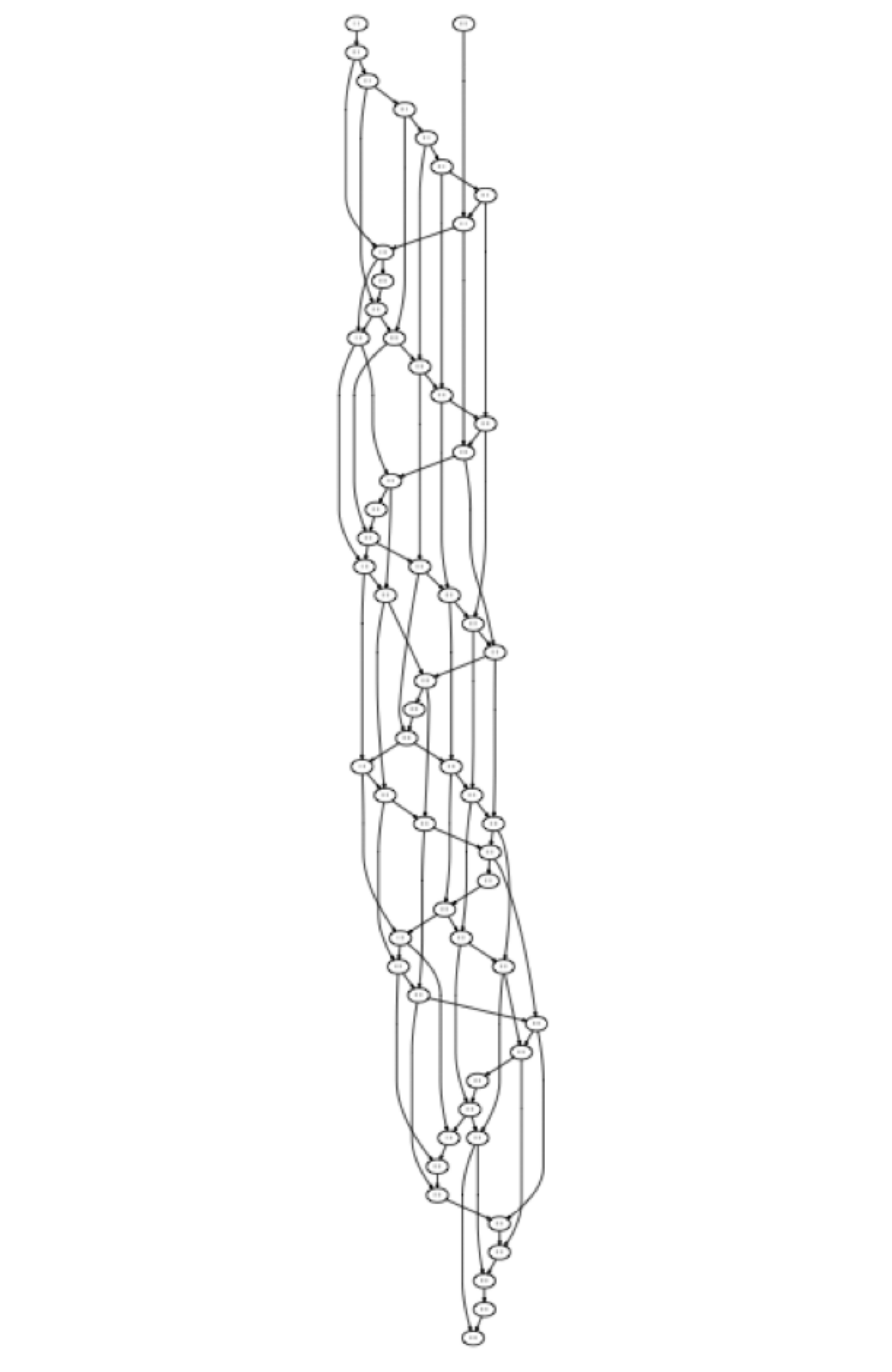}}
\hspace{20mm}
\subfigure[Non-mutual interaction.]{
\includegraphics[width=0.4\textwidth]{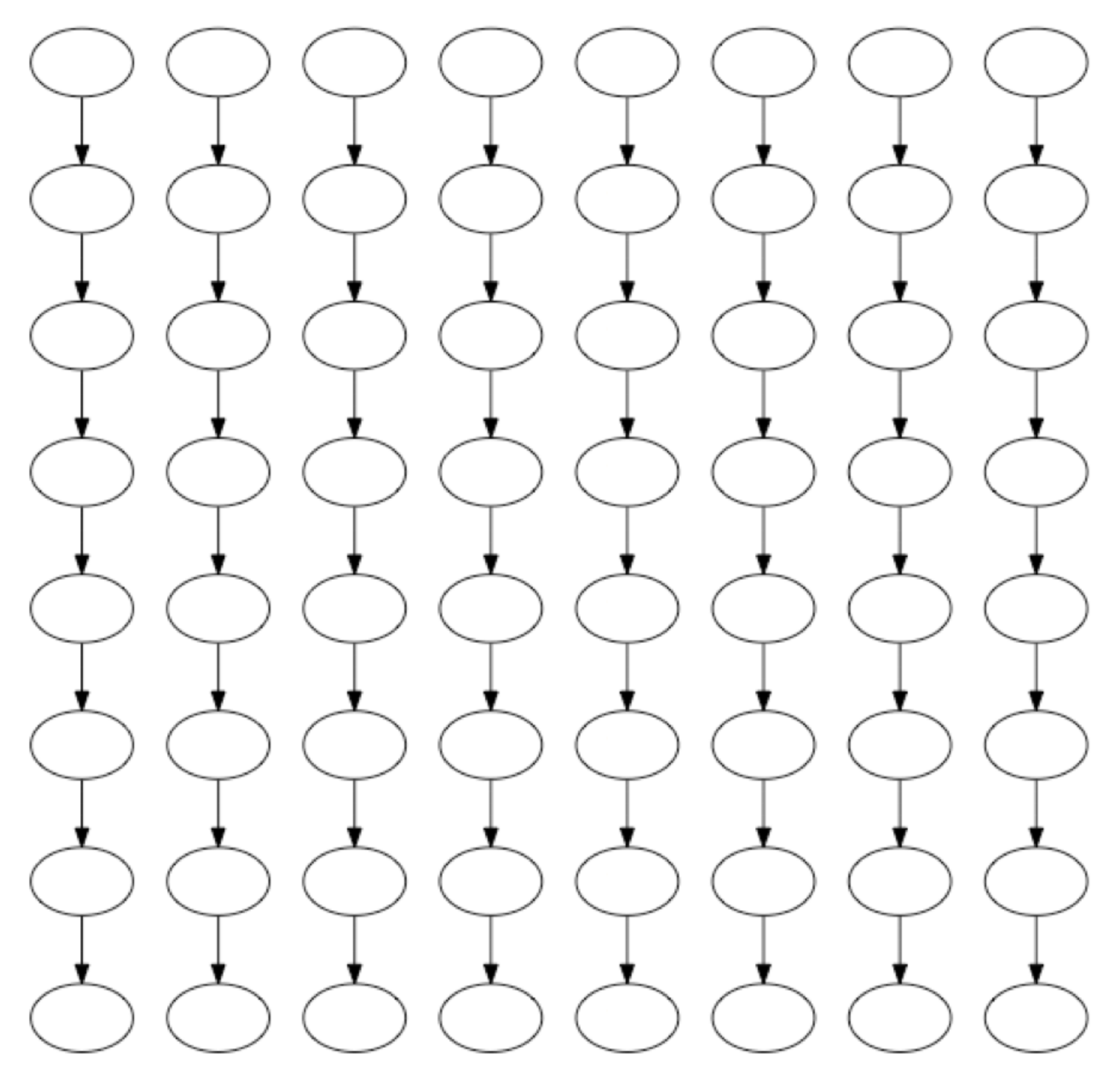}}
\caption{DAG of FMM traversal using mutual and non-mutual interactions.}
\label{fig:dag}
\end{figure}

In the FMM it is possible to calculate the \textit{mutual} interaction of cells, reducing the computation by about half. This is possible because all translation operators are a function of distance, and flipping the sign of the odd terms allows it to be reused when the target and source are interchanged. However, mutual interaction has a more restrictive data dependency, and Non-mutual interaction is more suitable for many cores due to the high degree of parallelism it provides. Figure~\ref{fig:dag} shows the actual DAG
of the FMM traversal using mutual or non-mutual particle interactions. The nodes represents the coarse-grained tasks and the edges the data dependencies between tasks. While the DAG for mutual interaction is very scattered, the DAG for non-mutual interaction clearly exposes more parallelism, which can be exploited as a low hanging
fruit by the QUARK dynamic scheduler. Furthermore, The \emph{LOCALITY} flag for the non-mutual
interaction in Figure~\ref{fig:fmm-mutual-interaction}
permits to ensure that the sequence of tasks within one branch is not interrupted by 
another thread during the execution, which may pollute the cache memory.\\

The main key for getting high performance is to tune the granularity of the \emph{super} task.
We have performed a thorough investigation of the 
optimal granularity, which will be presented in the following section.

\begin{figure}[tb]
\centering
\scriptsize
\begin{lstlisting}[language=C]
void Evaluator::interact(C_iter Ci, C_iter Cj, Quark *quark, bool mutual) {
  Quark_Task_Flags tflags = Quark_Task_Flags_Initializer;
  if( mutual ) {
    QUARK_Insert_Task(quark,interactQuark,&tflags,
                      sizeof(Cell),&*Ci,INOUT,
                      sizeof(Cell),&*Cj,INOUT,
                      sizeof(bool),&mutual,VALUE,
                      0);
  } else {
    QUARK_Insert_Task(quark,interactQuark,&tflags,
                      sizeof(Cell),&*Ci,OUTPUT | LOCALITY,
                      sizeof(Cell),&*Cj,NODEP,
                      sizeof(bool),&mutual,VALUE,
                      0);
  }
}
\end{lstlisting}
\caption{Example of inserting and executing a task in the scheduler
when computing cell interactions.
  The \emph{interact} routine inserts a task into the scheduler,
  passing it the sizes and pointers of arguments and their usages
  (INPUT, OUTPUT, INOUT, NODEP, VALUE) along with the LOCALITY flag
  for optimizing data reuse.}
\label{fig:fmm-mutual-interaction}
\end{figure}

\section{Experimental Results}
\label{sec:results}
In this section, we will present the results of a parametric scalability study of our data-driven FMM approach. All experiments are performed on a single node with a quad-socket quad-core \texttt{Intel\textregistered \ Xeon\textregistered \ E7340} processor at 2.40 GHz with 16 physical cores. The code is written in C++ and the compiler is \texttt{gcc version 4.1.2 20080704}, and OS is \texttt{x86\_64 Red Hat 4.1.2-50}. All codes were compiled with ``\texttt{g++ -O3 -fopenmp -ffast-math -funroll-loops -fforce-addr}".

We have mentioned in Section \ref{sec:impl} that a pair of subtrees are assigned to QUARK as tasks, instead of a pair of cells. Depending on how large these subtrees are, QUARK will have to balance either a large number of small tasks or a small number of large tasks through work stealing techniques. The idea here is to optimize the overhead of dynamic task scheduling, by adjusting the size of the subtrees that are passed to QUARK. In our \textit{dual tree traversal} approach, it is trivial to adjust the size of the subtrees that are passed to QUARK by simply stopping the breadth first (queue based) dual tree traversal, and sending all the pairs of subtrees in the queue to QUARK as individual tasks. We note that the use of a breadth first traversal is essential to this approach, since it creates tasks of similar size in the queue, whereas a depth first (stack based) traversal will result in largely varying task sizes in the stack, which would be much more difficult to balance. Furthermore, a breadth first traversal will create tasks with much higher data-parallelism. Also, the size of the queue can get much larger (tens of thousands) than the size of the stack in a depth first traversal (tens).

\begin{figure}[t!]
\subfigure[$p=6, N=10^5$]{
\includegraphics[width=0.45\textwidth]{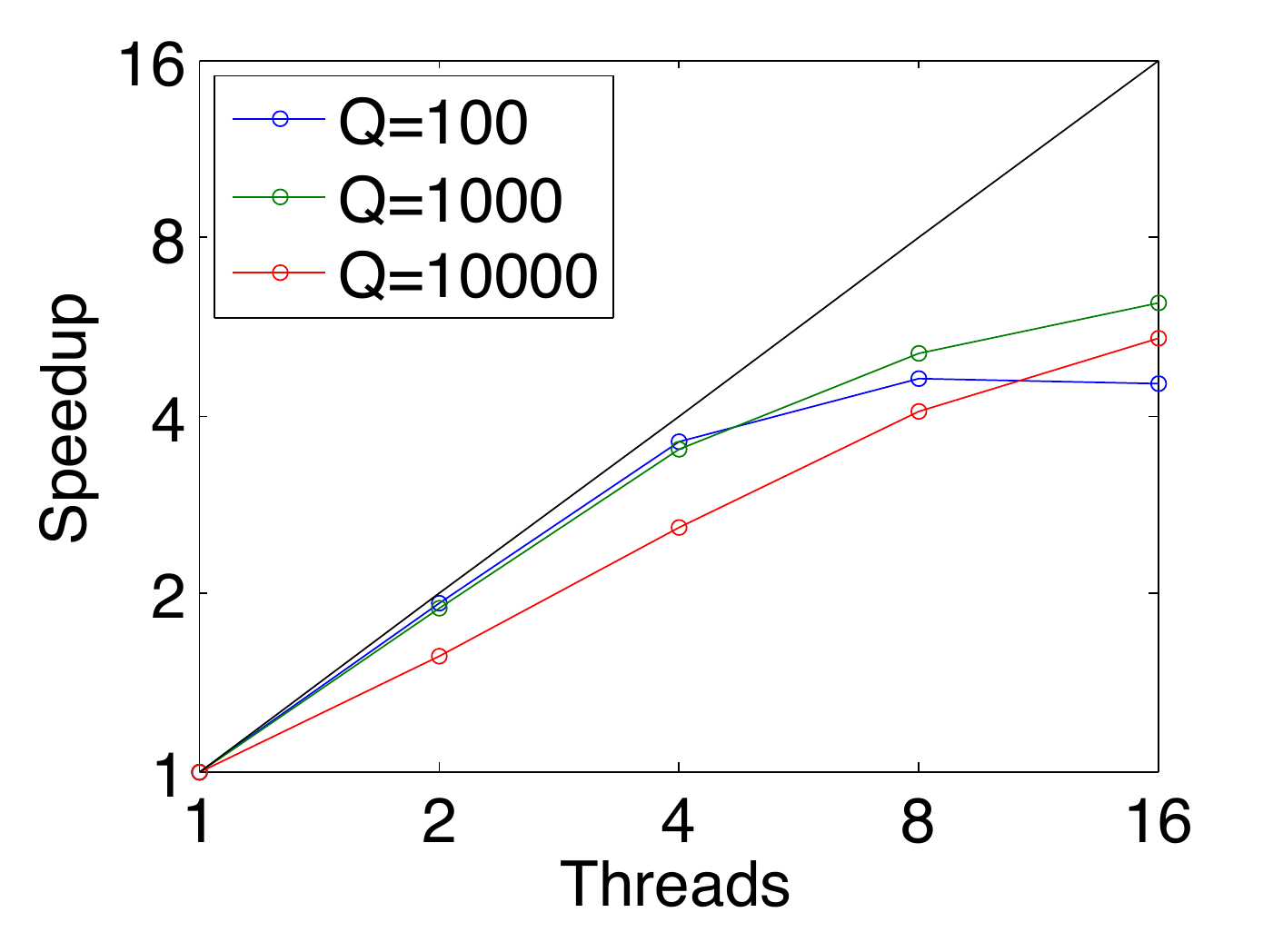}\label{fig:p6n5}}
\subfigure[$p=9, N=10^5$]{
\includegraphics[width=0.45\textwidth]{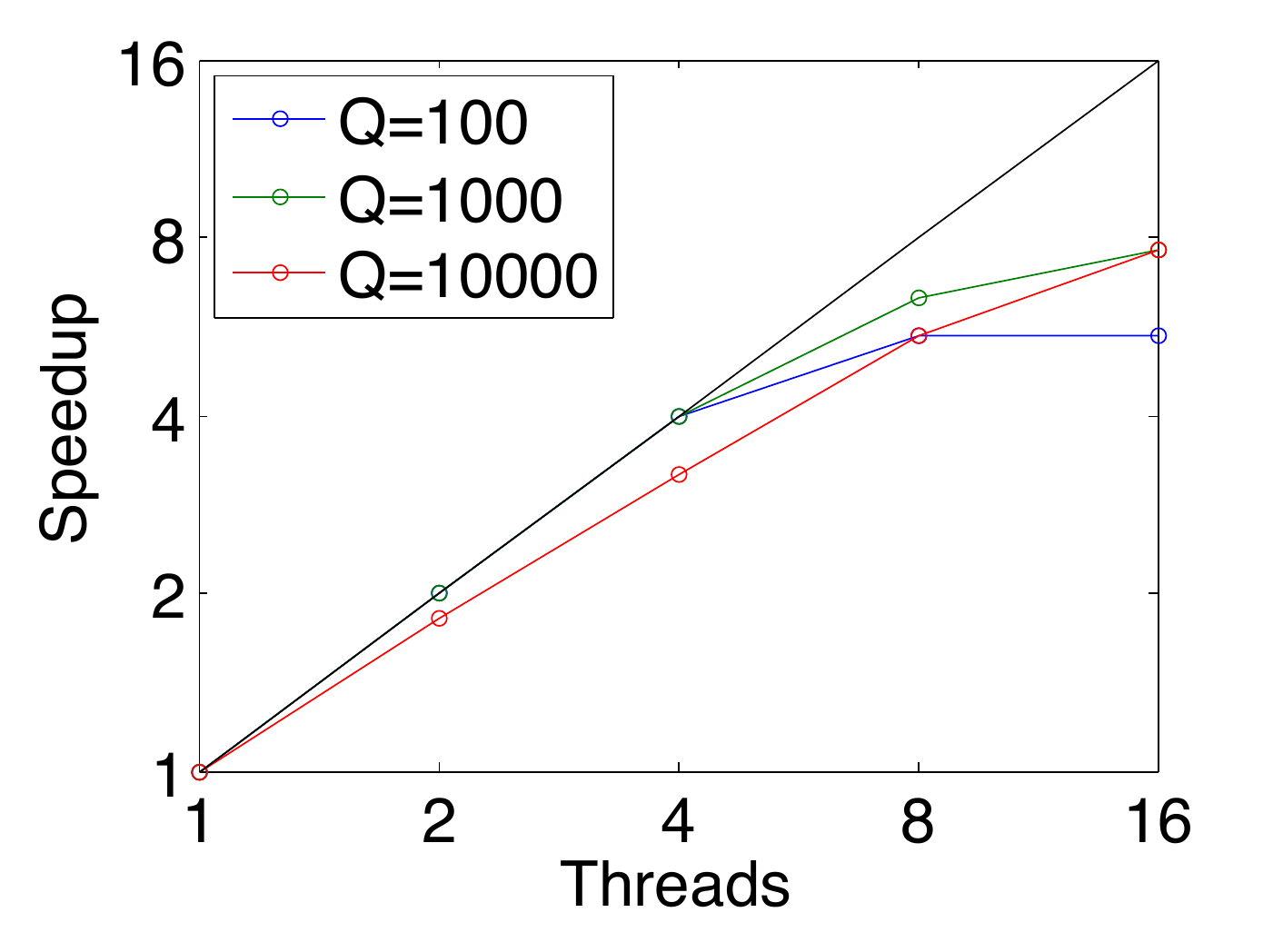}\label{fig:p9n5}}
\subfigure[$p=6, N=10^6$]{
\includegraphics[width=0.45\textwidth]{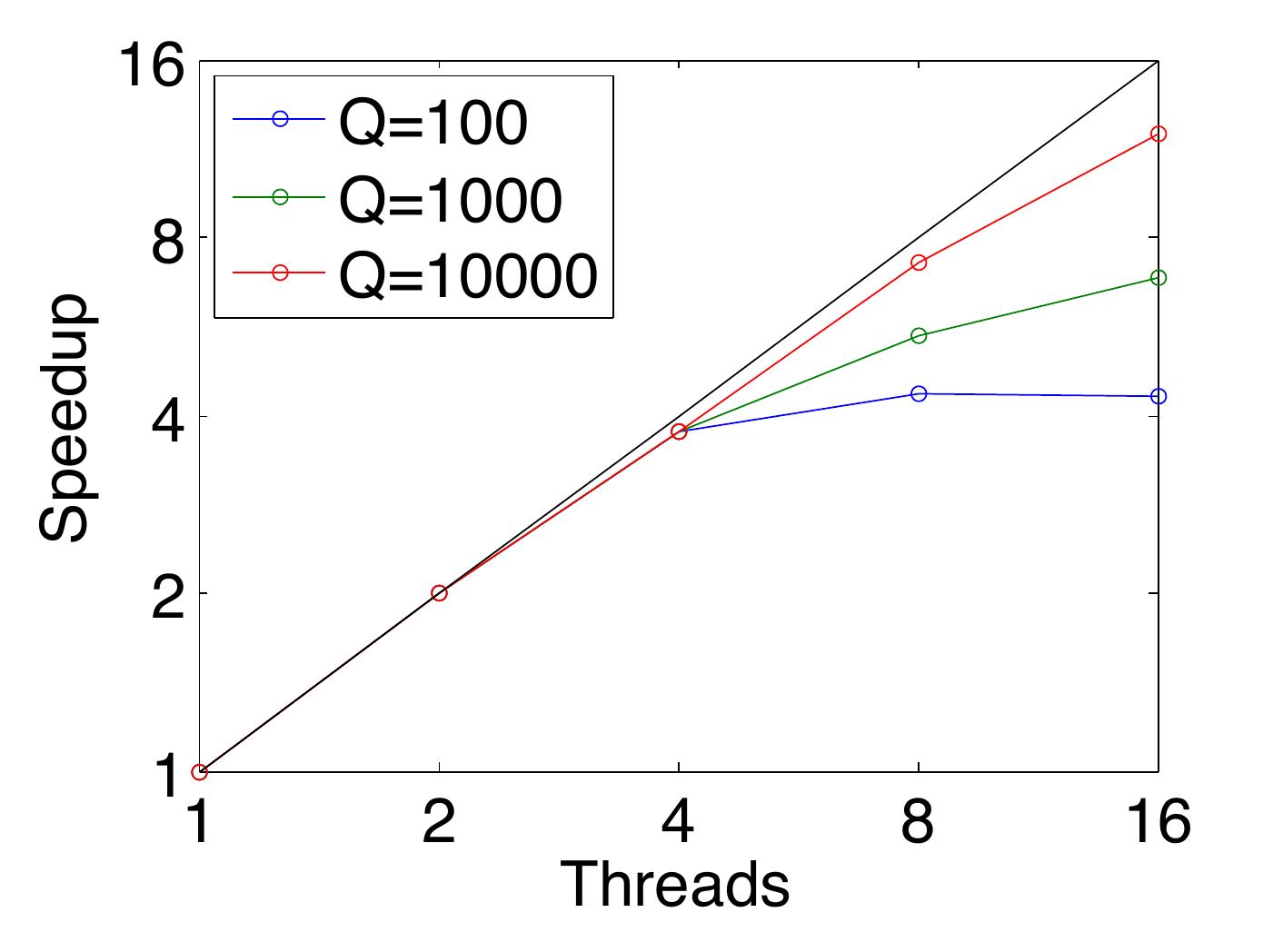}\label{fig:p6n6}}
\subfigure[$p=9, N=10^6$]{
\includegraphics[width=0.45\textwidth]{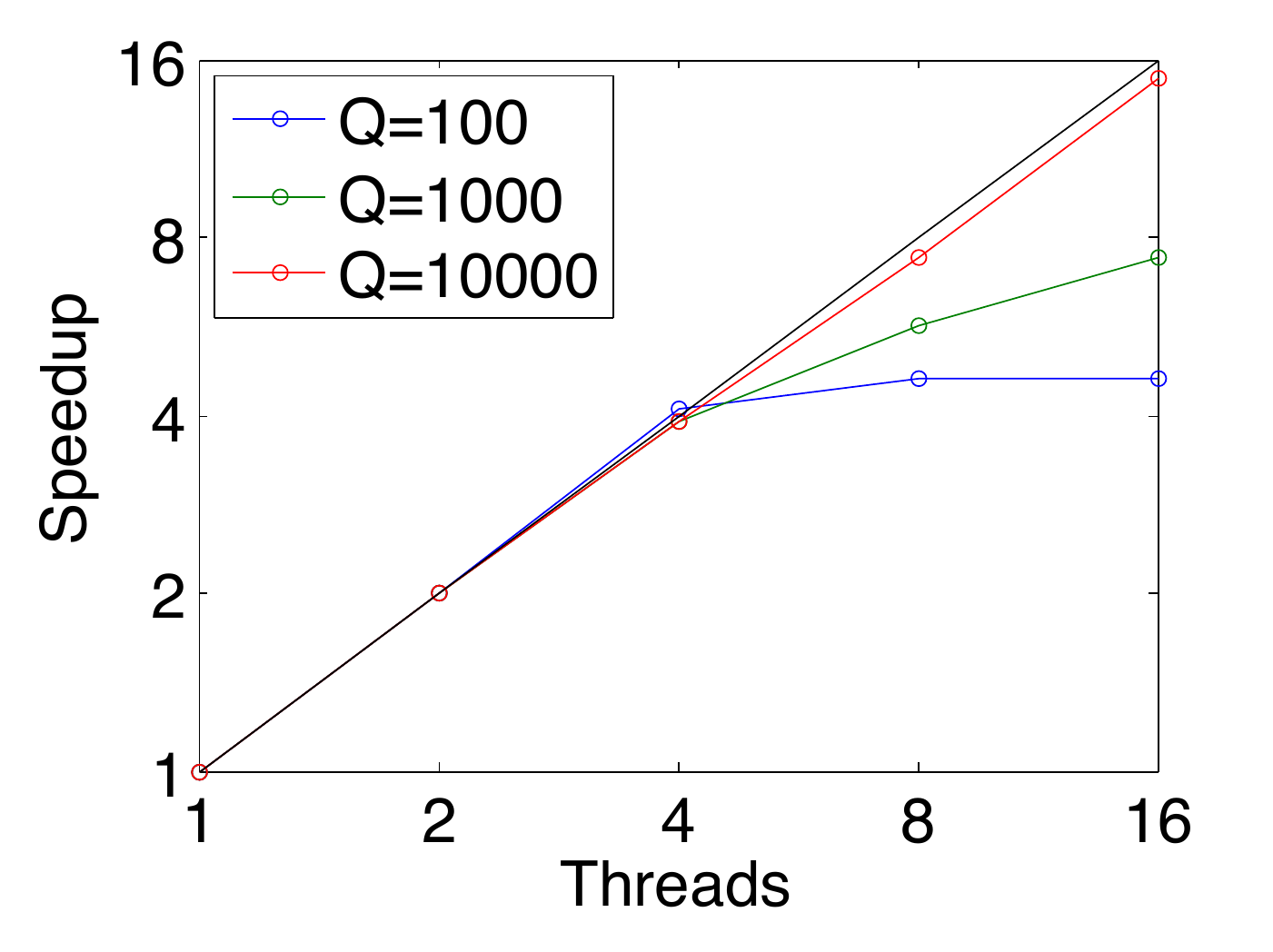}\label{fig:p9n6}}
\subfigure[$p=6, N=10^7$]{
\includegraphics[width=0.45\textwidth]{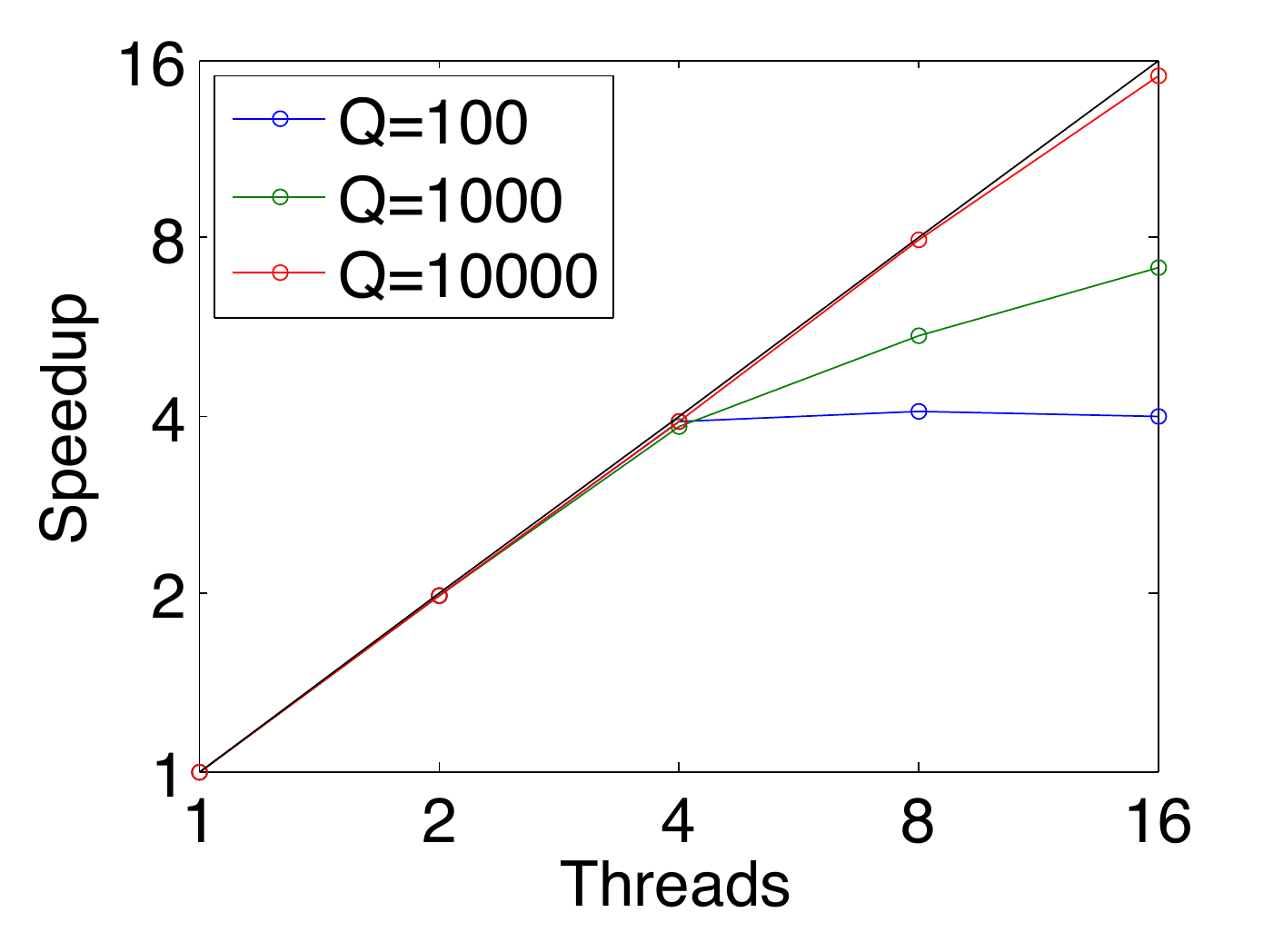}\label{fig:p6n7}}
\hspace{12mm}
\subfigure[$p=9, N=10^7$]{
\includegraphics[width=0.45\textwidth]{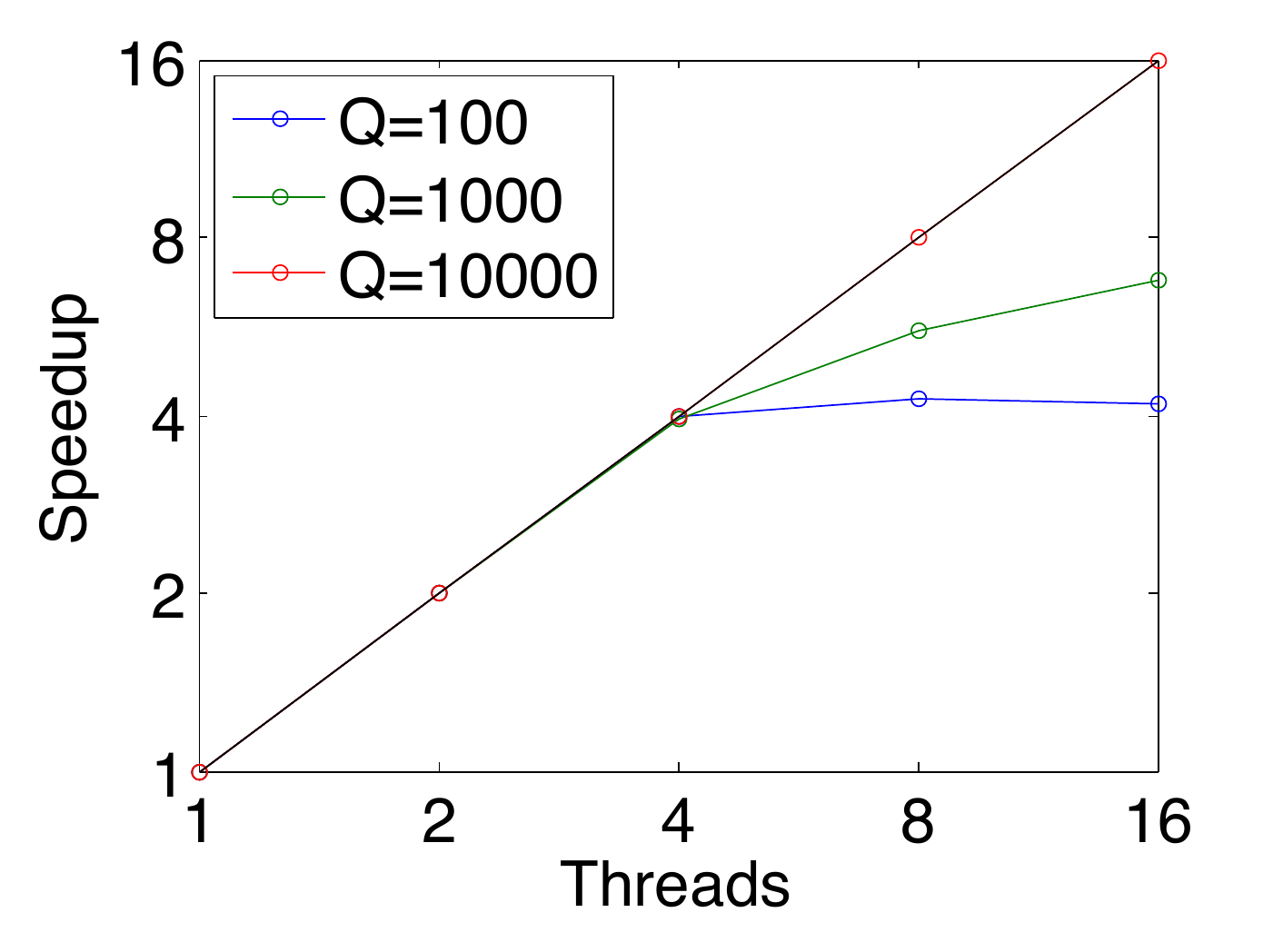}\label{fig:p9n7}}
\caption{Strong scaling of FMM on up to 16 threads. $p$ is the expansion order, $N$ is the number of particles, and $Q$ is the queue size which is dispatched to QUARK.}
\label{fig:scaling}
\end{figure}

The results of strong scaling tests for the FMM are shown in Fig. \ref{fig:scaling}, where $Q$ is the size of the queue, $p$ is the order of expansion, and $N$ is the number of particles. During the dual tree traversal the size of the queue $Q$ starts from $1$ (pair of root cells) and increases rapidly. We simply put a conditional statement in the tree traversal to ship all the queued tasks to QUARK when the queue size reaches a certain threshold, and this threshold corresponds to the value $Q$. The particles are randomly distributed in a unit cube and the resulting tree is well balanced.

By comparing Figures \ref{fig:p6n5}, \ref{fig:p6n6}, and \ref{fig:p6n7} we see that increasing the problem size $N$ automatically gives us better strong scalability, but only if $Q$ is sufficiently large. When the task size is small, the overhead of dynamically scheduling $10,000$ tasks becomes a significant burden as shown in Figures \ref{fig:p6n5} and \ref{fig:p9n5}. However, as the task size increases proportional to $N$, the constant overhead of scheduling $10,000$ tasks can be amortized as in Figures \ref{fig:p6n7} and \ref{fig:p9n7}.
From the fact that $Q=10,000$ is too large for $N=10^5$ but not for $N=10^6$, we may assume that a ratio in the order of $N/Q>100$ is necessary to amortize the cost of the dynamic scheduling. Furthermore, increasing the order of expansion $p$ increases the operational intensity of the M2L kernel and results in slightly better scalability, which can be seen by comparing the left and right column of Figure \ref{fig:scaling}. In Figure \ref{fig:p9n7} we achieve 100 \% parallel efficiency on 16 cores, using $p=9$, $N=10^7$, and $Q=10,000$.

\begin{figure}[t]
\centering
\includegraphics[width=0.8\textwidth]{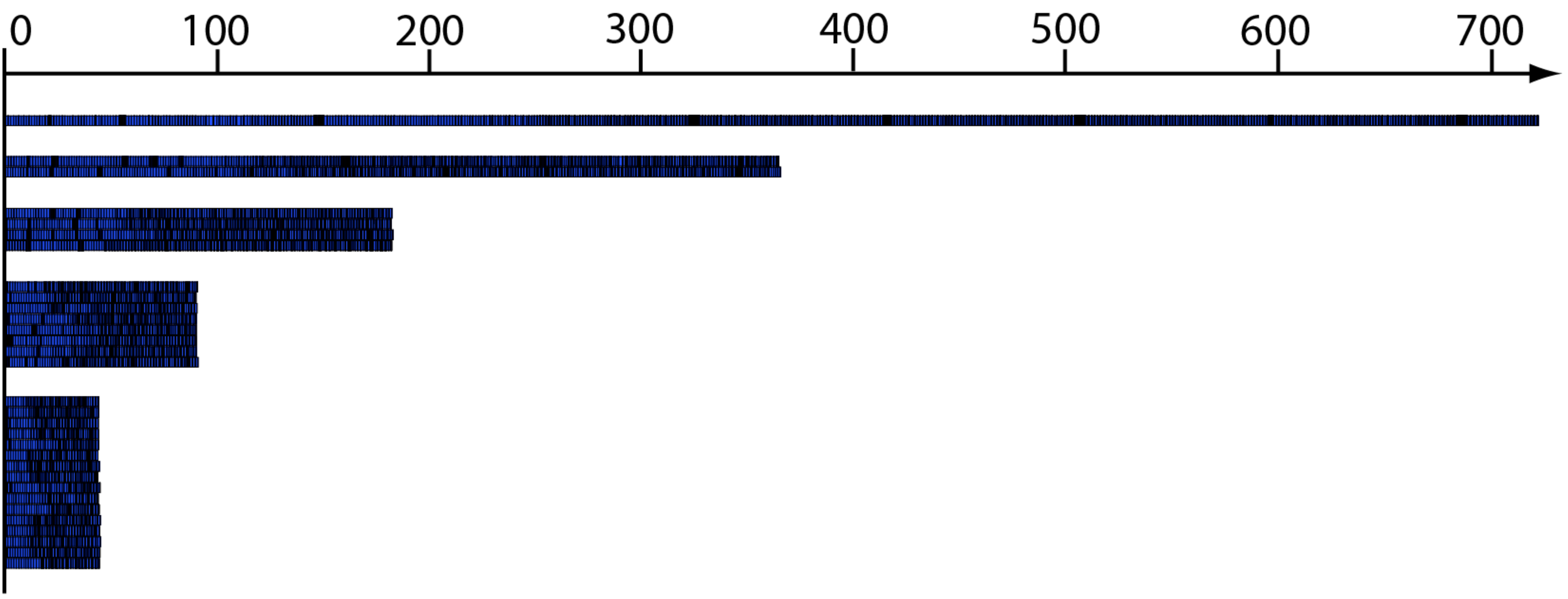}
\caption{Execution trace of the data-driven FMM using 1, 2, 4, 8, and 16 threads.}
\label{fig:trace}
\end{figure}

The load-imbalance in the thread execution for $p=9$, $N=10^7$, and $Q=10,000$ is shown in Figure \ref{fig:trace}. The tick marks represent the total runtime of the FMM in seconds. The runtime is 720, 360, 180, 90, and 45 seconds on 1, 2, 4, 8, and 16 threads, respectively. Each blue segment represents an individual task that is dynamically scheduled by QUARK. These tasks are composed of dual tree traversals of the pair of subtrees and the corresponding M2L/P2P kernels that arise from these traversals. Since the searching of \textit{well separated} cell pairs is imbedded inside each task unit and scheduled dynamically, having an irregular tree structure should not degrade the quality of the load balance.

\section{Conclusion and Future Work}
\label{sec:conclusion}
This paper describes a data-driven execution of FMMs based on the dynamic runtime
system QUARK. After carefully tuning the granularity of the subtrees, our implementation
achieves a linear speedup on a quad-socket quad-core Intel Xeon (16 cores total).
QUARK permits to achieve not only high performance but also high productivity
in terms of parallel implementation. The end user can therefore focus 
and spend time improving his core numerical kernels and the burden to get
parallel performance is rather shifted on the runtime.

The authors plan to extend this work with the StarPU~\cite{Augonnet_2010_ccpe}
dynamic runtime, which schedules tasks on x86 as well as on hardware
accelerators (GPUs). Although it has a somewhat similar API than QUARK, 
the user has still to develop the appropriate kernels for the GPU
and can let it up to the runtime to decide on which available resource, x86 or GPU,
the task can be executed on. Moreover, StarPU provides a reduction operation,
which could further improve our current implementation by
adding another dimension of parallelism during the execution of the 
successive tasks seen in the non-mutual DAG representation. Finally,
the authors will eventually tackle the distributed memory environment  
using the DAGuE~\cite{DAGuE} framework to perform load balancing across time steps.

\bibliographystyle{abbrv}

\end{document}